\input amstex
\input epsf
\magnification=\magstep1 
\baselineskip=13pt
\documentstyle{amsppt}
\vsize=8.7truein \CenteredTagsOnSplits \NoRunningHeads
\def\ii{\bold{i}}

\def\rk{\operatorname{rank}}
\def\spa{\operatorname{span}}
\def\dist{\operatorname{dist}}

\def\Pr{\bold {P}}
\topmatter
 
\title Computing the theta function   \endtitle 
\author Alexander Barvinok  \endauthor
\address Department of Mathematics, University of Michigan, Ann Arbor,
MI 48109-1043, USA \endaddress
\email barvinok$\@$umich.edu \endemail
\date June 4,  2023 \enddate
\thanks  This research was partially supported by NSF Grant DMS 1855428. 
\endthanks 
\abstract Let $f: {\Bbb R}^n \longrightarrow {\Bbb R}$ be a positive definite quadratic form and let $y \in {\Bbb R}^n$ be a point. We present a fully polynomial randomized approximation scheme (FPRAS) for computing $\sum_{x \in {\Bbb Z}^n} e^{-f(x)}$, provided the eigenvalues of $f$ lie in the interval roughly between $s$ and $e^{s}$ and for computing $\sum_{x \in {\Bbb Z}^n} e^{-f(x-y)}$, provided the eigenvalues of $f$ lie in the interval roughly between $e^{-s}$ and $s^{-1}$ for some $s \geq 3$. To compute the first sum, we represent it as the integral of an explicit log-concave function on ${\Bbb R}^n$, and to compute the second sum, we use the reciprocity relation for theta functions. We then apply our results to test the existence of many short integer vectors in a given subspace $L \subset {\Bbb R}^n$, to estimate the distance from a given point to a lattice, and to sample a random lattice point from the discrete Gaussian distribution.
 \endabstract
\endtopmatter
\document

\head 1. Introduction \endhead

\subhead (1.1) The theta function \endsubhead Let $f: {\Bbb R}^n \longrightarrow {\Bbb R}_+$ be a positive definite quadratic form, so 
$$f(x) = \langle B x, x \rangle \quad \text{for} \quad x \in {\Bbb R}^n,$$
where $B$ is an $n \times n$ positive definite matrix and $\langle \cdot, \cdot \rangle$ is the standard scalar product in ${\Bbb R}^n$.
We consider the problem of efficient computing (approximating) the sum 
$$\Theta(B)=\sum_{x \in {\Bbb Z}^n} e^{-f(x)} = \sum_{x \in {\Bbb Z}^n} e^{-\langle Bx, x \rangle}, \tag1.1.1$$
where ${\Bbb Z}^n \subset {\Bbb R}^n$ is the standard integer lattice. More generally, for a given point $y \in {\Bbb R}^n$, we want to efficiently compute (approximate) the sum
$$\Theta(B, y)=\sum_{x \in {\Bbb Z}^n} e^{-f(x-y)} = \sum_{x \in {\Bbb Z}^n} e^{- \langle B(x-y), x-y \rangle}. \tag1.1.2$$
Together with (1.1.1) and (1.1.2), we also compute the sum
$$\sum_{x \in {\Bbb Z}^n} \exp\left\{ - \langle Bx, x \rangle +  \ii  \langle b, x \rangle \right\}, \tag1.1.3$$
where $b \in {\Bbb R}^n$ and $\ii^2=-1$.

Of course, the sums (1.1.1) -- (1.1.3) are examples of the (multivariate) theta function, an immensely popular object, 
see, for example, \cite{M07a}, \cite{M07b} and \cite{M07c}. Theta functions satisfy the {\it reciprocity relation}
$$\split & \sum_{x \in {\Bbb Z}^n} \exp\left\{- \pi \langle B(x -y), x-y\rangle \right\} \\=&{1 \over\sqrt{\det B}} \sum_{x \in {\Bbb Z}^n} 
\exp\left\{- \pi \langle B^{-1} x,  x \rangle + 2 \pi \ii \langle  x, y \rangle\right\}, \endsplit \tag1.1.4$$
see, for example, \cite{BL61}.

One motivation to study  (1.1.1)--(1.1.3) from the computational point of view comes from connections with algorithmic problems on lattices,
such as approximating the length of a shortest non-zero vector in the lattice and estimating the distance from a given point to a given lattice, see \cite{Sc87}, \cite{G+93}, \cite{Ba93}, \cite{Aj96}, \cite{A+01}, \cite{MG02}, \cite{D+03}, \cite{AR05}, \cite{Kh05}, \cite{MR07}, \cite{A+15}, \cite{M+21}, as well as lattice-based cryptography, see \cite{MG02}, \cite{MR07}, \cite{G+08}, \cite{MR09}, \cite{Pe10}.

\subhead (1.2) Lattices \endsubhead A {\it lattice} $\Lambda \subset {\Bbb R}^n$ is a discrete additive subgroup which spans ${\Bbb R}^n$. Equivalently, 
$\Lambda$ is the set of all integer linear combinations of some linearly independent vectors $u_1, \ldots, u_n$, called a {\it basis} of $\Lambda$,
$$\Lambda=\left\{ \sum_{i=1}^n \xi_i u_i: \quad \xi_i \in {\Bbb Z} \quad \text{for} \quad i=1, \ldots, n \right\}.$$
We say that $\rk \Lambda=n$.

For $n >1$, the same lattice $\Lambda$ has many different bases, and some of those bases are more convenient to work with than others, see, for example, 
\cite{G+93} and \cite{MG02}. Given a vector  $u \in \Lambda$, $u=\xi_1 u_1 + \ldots + \xi_n u_n$, we have 
$$\|u\|^2 = \langle Bx, x \rangle \quad \text{where} \quad x=\left(\xi_1, \ldots, \xi_n\right)$$
and $B$ is the Gram matrix of the vectors $u_1, \ldots, u_n$, so that 
$$B=\left(\beta_{ij}\right) \quad \text{where} \quad \beta_{ij}=\langle u_i, u_j \rangle.$$
Similarly, if $v \in {\Bbb R}^n$ is an arbitrary point, $v=\eta_1 u_1 + \ldots  + \eta_n u_n$, then 
$$\|u-v\|^2 = \langle B(x-y), x-y \rangle \quad \text{for} \quad y=\left(\eta_1, \ldots, \eta_n\right).$$
Consequently, the theta functions (1.1.1) and (1.1.2) are written as
$$\Theta(B)=\sum_{u \in \Lambda} e^{-\|u\|^2} \qquad \text{and} \qquad \Theta(B, y)=\sum_{u \in \Lambda} e^{-\|u-v\|^2}. \tag1.2.1$$
We see from (1.2.1) that the theta functions do not depend on the choice of a basis of the lattice: choosing a different basis corresponds to replacing 
the Gram matrix $B$ with a Gram matrix of the form $A^T B A$, where $A \in GL(n, {\Bbb Z})$ is an integer matrix such that $\det A =\pm 1$. It follows that the value of $\det B$ does not depend on the choice of a basis. 
The number $\sqrt{\det B}$ is called the {\it determinant} of $\Lambda$ and denoted $\det \Lambda$, see, for example, Chapter I of \cite{Ca97}.

The following two optimization problems have attracted a lot of attention due to their importance for optimization and cryptography. One is finding (or approximating) the minimum length of a non-zero vector 
from a given lattice,
$$\lambda(\Lambda)=\min_{u \in \Lambda \setminus \{0\}} \|u\|$$
and the other is finding (or approximating) the distance from a given point $v \in {\Bbb R}^n$ to a given lattice,
$$\dist(v, \Lambda) = \min_{u \in \Lambda} \|u-v\|,$$
see \cite{Sc87}, \cite{G+93}, \cite{Ba93}, \cite{Aj96}, \cite{A+01}, \cite{MG02}, \cite{D+03}, \cite{AR05}, \cite{Kh05},  \cite{A+15}. We  assume that $\Lambda$ is defined by its basis.
In a breakthrough paper \cite{Ba93}, Banaszczyk used theta functions to obtain structural results (known as ``transference theorems") for $\lambda(\Lambda)$ and a host of related quantities (successive minima, covering radius, etc.) 
Using results of \cite{Ba93}, Aharonov and Regev \cite{AR05} showed that the problems of approximating $\lambda(\Lambda)$ and 
$\dist(v, \Lambda)$ within a factor $O(\sqrt{n})$ lie in NP $\cap$ co-NP. This is in contrast to the fact that the existing polynomial time algorithms are guaranteed to approximate the desired quantities roughly within a $2^{O(n)}$ factor, more precisely within a factor of $2^{O(n (\log \log n)^2/\log n)}$ in deterministic polynomial time \cite{Sc87} and within a factor 
$2^{O(n \log \log n/\log n)}$ in randomized polynomial time \cite{A+01}. Computing $\lambda(\Lambda)$ exactly is NP-hard, and approximating $\lambda(\Lambda)$ within a factor of $2^{(\log n)^{{1 \over 2}-\epsilon}}$ is hard modulo some plausible computational complexity assumptions \cite{Kh05},
while approximating $\dist(v, \Lambda)$ within a factor of $n^{c/\log \log n}$ is NP-hard for some absolute constant $c >0$ \cite{D+03}.

Given a lattice $\Lambda \subset {\Bbb R}^n$ and a point $v \in {\Bbb R}^n$, one can define a probability measure on $\Lambda$, called 
the {\it discrete Gaussian distribution}, where the probability of $u \in \Lambda$ is proportional to $e^{-\|u-v\|^2}$,
$$\Pr(u) \sim e^{-\|u-v\|^2} \quad \text{for all} \quad u \in \Lambda. \tag1.2.2$$
Efficient approximate sampling from the distribution (1.2.2) has attracted a lot of attention, in connection with optimization and cryptography, see
\cite{G+08}, \cite{MR07}, \cite{MR09}, \cite{Pe10}, \cite{A+15}, \cite{RS17}.

\head 2. Results \endhead 

\subhead (2.1) Approximating the theta function \endsubhead 
In what follows, we write $A \preceq B$ for $n \times n$ real symmetric matrices $A$ and $B$ if $B-A$ is a positive semidefinite matrix. We denote by $I$ the $n \times n$ identity matrix. Our main result is a fully polynomial randomized approximation scheme (FPRAS) for computing (1.1.1) and 
(1.1.3) provided 
$$sI \ \preceq \ B \ \preceq \ \left(s+ {e^s \over 4} \left(1-e^{-s}\right)^2 \left(1-e^{-2s}\right)\right) I \quad \text{for some} \quad s \geq 1. \tag2.1.1$$
Thus we present a randomized algorithm that for any $B$ satisfying (2.1.1) and for any $\epsilon > 0$ approximates the value of $\Theta(B)$ and that of (1.1.3) within relative error $\epsilon$ in time polynomial in $n$, $\epsilon^{-1}$ and $s$. 
It turns out that when (2.1.1) is satisfied, we can write (1.1.1) and (1.1.3)  as an integral of some explicit log-concave function $G: {\Bbb R}^n \longrightarrow {\Bbb R}_+$ and hence we can use any of the efficient algorithms for integrating log-concave functions as a blackbox \cite{AK91}, \cite{F+94}, \cite{FK99}, \cite{LV07}.  From (2.1.1) we obtain an easier to parse condition 
$$s I \ \preceq \ B \ \preceq \ \left(s + {e^s \over 5}\right) I \quad \text{for} \quad s \geq 3, \tag2.1.2$$
which is sufficient for $\Theta(B)$ and, more generally, for (1.1.3) to be efficiently computable. We describe the algorithm is Section 3 and prove the main structural result (Theorem 3.1) underlying the algorithm in Section 4.

From the reciprocity relation (1.1.4) it immediately follows that there is an FPRAS for $\Theta(B, y)$ provided 
$$\aligned &\pi^2 \left(s+ {e^s \over 4} \left(1-e^{-s}\right)^2 \left(1-e^{-2s}\right)\right)^{-1} I  \ \preceq \ B \ \preceq \ \pi^2 s^{-1} I \\ &\qquad \text{for some} \quad s \geq 1. \endaligned
\tag2.1.3$$
That is, there is a randomized algorithm that for any $B$ satisfying (2.1.3), for any $y \in {\Bbb R}^n$ and any $0 < \epsilon < 1$ approximates the value of $\Theta(B, y)$ within relative error $\epsilon$ in time polynomial in $n$, $\epsilon^{-1}$ and $s$.
An easier to parse sufficient condition is 
$$\pi^2 \left( s +{e^s \over 5}\right)^{-1} I \ \preceq \ B \ \preceq \ \left(\pi^2s^{-1}\right) I \quad \text{for} \quad s \geq 3. \tag2.1.4 $$

\remark{(2.1.5) The smooth range}

Let us fix $\gamma > 1$ and let $s=\gamma \ln n$.
It is not hard to check that if $sI \preceq B$ then the value of  $\Theta(B)$, and, more generally, of (1.1.3) is $1+ O\left(n^{1-\gamma}\right)$, since only $x=0$ contributes significantly to the sum. Furthermore, a straightforward algorithm approximates $\Theta(B)$ and (1.1.3) within relative error $\epsilon$ in time polynomial in $n$ and $\epsilon^{-1}$, provided $n$ is sufficiently large, $n \geq n_0(\gamma)$. For the sake of completeness, we present the algorithm along with some technical estimates in Section 8.

Applying the reciprocity relation (1.1.4), we have 
$$\Theta(B,y) = {\pi^{n/2} \over \sqrt{\det B}} \left(1+ O\left(n^{1-\gamma}\right)\right) \quad \text{provided} \quad B \preceq \left({\pi^2 \over \gamma \ln n}\right) I
\quad \text{for} \quad \gamma >1.$$
Furthermore, as long as $\gamma >1$ is fixed, for any $\epsilon >0$ the value of $\Theta(B, y)$ can be approximated within relative error $\epsilon$ in time polynomial in $n$ and $\epsilon^{-1}$.
Hence if $B$ is sufficiently small in the ``$\preceq$" order, the discrete sum (1.1.2) is well-approximated by the integral 
$$\int_{{\Bbb R}^n}
\exp\left\{- \langle B(x-y), x-y \rangle\right\} \ dx = {\pi^{n/2} \over \sqrt{\det B}}.$$ This phenomenon is described by the {\it smoothing parameter} of a lattice introduced in \cite{MR07}. Our constraints (2.1.1) and (2.1.3) correspond to the ``non-smooth" range when $s \leq \gamma \ln n$ for some fixed 
$0 < \gamma < 1$. Apart from some straightforward situations (for example, when the matrix $B$ is diagonal), the condition (2.1.3) appears to be the first one when $\Theta(B, y)$ can be efficiently approximated in a non-smooth, that is genuinely discrete, case.  
\endremark

\subhead (2.2) Integer points in a subspace \endsubhead Let $A$ be an $m \times n$ integer matrix of $\rk A=m < n$ and let 
$$\Lambda = \left\{ x \in {\Bbb Z}^n: \ Ax=0 \right\}. \tag2.2.1$$
Then $\Lambda$ is a lattice in the ambient space $\spa(\Lambda)=\ker A$. We remark that even when $m=1$, the class of such lattices (2.2.1) is quite rich: it is shown in \cite{S+11} that any lattice $\Lambda'$ of rank $n$ can be arbitrarily closely approximated by a proper scaling $\alpha \Lambda$ of a lattice $\Lambda$ that is a hyperplane section of ${\Bbb Z}^{n+1}$.

For $s >0$, we consider the theta function
$$\Theta_{\Lambda}(s)=\sum_{u \in \Lambda} e^{-s \|u\|^2}. $$
We denote by $\|A\|_{\text{op}}$ the operator norm of $A$, that is the largest singular value of $A$.
Let us fix $\delta >0$. In what follows, we consider asymptotics as $n$ grows.

 In Section 5, we show that if $\|A\|_{\text{op}} =o\left(n^{\delta}\right)$, then for 
$$s=\left({1 \over 2} +\delta\right) \ln n$$
and any $\epsilon > 0$, the value of $\Theta_{\Lambda}(s)$ can be approximated within relative error $\epsilon+o(1)$ in randomized polynomial time. This is based on the observation that $\Theta_{\Lambda}(s)$ is approximated within an additive error $o(1)$ by the function $\Theta(B)$ of (1.1.1), where $B$ is an $n\times n$ matrix with the eigenvectors in $\ker A$ with eigenvalue $s$ and in $(\ker A)^{\bot} = \operatorname{im} A^T$ with eigenvalue $s + e^s/5$ so that $B$ satisfies 
(2.1.2) when $s \geq 3$.

Note that as long as $\delta < 1/2$, we are in a ``non-smooth" range, cf. Section 2.1.5.

This result is then applied to testing the existence of short non-zero vectors in $\Lambda$. We show that if 
$$\min_{u \in \Lambda \setminus\{0\}} \|u\| \ \gg \ n^{{1 \over 2}-\delta}$$
then $\Theta_{\Lambda}(s)=1 + o(1)$,
while $\Theta_{\Lambda}(s) \gg 1$, if $\Lambda$ contains many short vectors, which allows us to separate these two cases in randomized polynomial time. 

Using a different approach, in \cite{M+21}, the authors present a polynomial time algorithm to find a lattice vector closest to a given point, when $A$ is a totally unimodular matrix.

\subhead (2.3) Estimating the distance to the lattice \endsubhead In Section 6, we consider the problem of estimating the distance from a given point 
$v \in {\Bbb R}^n$ to a given lattice $\Lambda \subset {\Bbb R}^n$, provided ${\Bbb Z}^n \subset \Lambda$. Such lattices $\Lambda$ appear in a few natural ways. If $\Lambda_0 \subset {\Bbb Z}^n$ is a lattice with an integer basis, then the {\it dual} or {\it reciprocal} lattice $\Lambda=\Lambda_0^{\ast}$ defined by 
$$\Lambda_0^{\ast}=\left\{ u \in {\Bbb R}^n: \quad \langle u, w \rangle \in {\Bbb Z} \quad \text{for all} \quad w \in \Lambda_0 \right\}$$
 contains ${\Bbb Z}^n$. The {\it $q$-ary} lattices $\Lambda$ satisfying $\left(q{\Bbb Z}\right)^n \subset \Lambda \subset {\Bbb Z}^n$ for an integer $q>1$ play a prominent role in lattice-based cryptography, see \cite{Aj96}, \cite{MG02}, \cite{MR09}. Typically, they are defined as the sets of solutions to systems of integer linear equations $\mod q$. Clearly, if $\Lambda$ is a $q$-ary lattice then the lattice $q^{-1} \Lambda$ contains ${\Bbb Z}^n$. 

For a lattice $\Lambda \subset {\Bbb R}^n$ and $\tau >0$, we define 
$$\Theta_{\Lambda}(\tau, v) =\sum_{u \in \Lambda} \exp\left\{ - \tau \|u-v\|^2 \right\}.$$ In particular, if $\Lambda ={\Bbb Z}^n$, then 
$\Theta_{{\Bbb Z}^n}(\tau, v)=\Theta(\tau I, v)$ and $\Theta_{{\Bbb Z}^n}(\tau, 0)=\Theta(\tau I)$ in the notation of Section 1.1. 

In Section 6, we prove that if ${\Bbb Z}^n \subset \Lambda$ then for any  $0 < \tau \leq 1$, we have 
$$41 e^{-\pi^2/\tau} \dist^2(v, \Lambda) \ \geq \ \ln {\Theta(\tau I) \over \Theta_{\Lambda}(\tau, v)} \ \geq \ 13 e^{-\pi^2/\tau} \dist^2(v, \Lambda) + \ln \det \Lambda. \tag2.3.1$$
As $n$ grows, under some conditions the additive term of $\ln \det \Lambda$ becomes asymptotically negligible and (2.3.1) provides an approximation of $\dist(v, \Lambda)$ within a constant factor of 
$\sqrt{41/13} \approx 1.8$, computable in randomized polynomial time. We provide an example to that effect in Section 6.
 
\subhead (2.4) Sampling from the discrete Gaussian distribution \endsubhead  Given a lattice $\Lambda \subset {\Bbb R}^n$ and a point $v \in {\Bbb R}^n$, 
we consider the discrete Gaussian probability distribution (1.2.2). Suppose that $\Lambda$ has a basis whose Gram matrix $B$ satisfies
$$\lambda I \ \preceq \ B \tag2.4.1$$
for some $\lambda >0$. Assume further that for any given $y \in {\Bbb R}^n$, the value of $\Theta(B, y)$ can be approximated in randomized polynomial time (for example, if $B$ satisfies (2.1.3)). We present an algorithm which for any given $0 < \epsilon < 1$ samples a random point $u \in \Lambda$ from a probability distribution $\mu$ which is $\epsilon$-close to (1.2.2) in the total variation distance, that is, 
$${1 \over 2} \sum_{u \in \Lambda} \left| \Pr(u)-\mu(u)\right| \ \leq \ \epsilon.$$
The complexity of the algorithm is polynomial in $n$, $\epsilon^{-1}$ and $\lambda^{-1}$.

It appears that previously polynomial time sampling algorithms, apart from some simple cases (such as when $B$ is a diagonal matrix), were known only in the smooth range, when the discrete Gaussian measure is well-approximated by its classical continuous version \cite{G+08}, \cite{Pe10}. Our algorithm follows the general logic of Peikert's algorithm \cite{Pe10}, except that we are able to extend it to non-smooth cases, since we are able to approximate the value of the theta function in those cases. Apart from that, the price we apparently have to pay is the dependence of the computational complexity on $\lambda$ in (2.4.1), which is absent in the smooth case.

We discuss the algorithm in Section 7.

\subhead (2.5) The plan of the paper \endsubhead Summarizing, the plan of the paper is as follows.

In Section 3, we present our main algorithm for approximating the theta functions (1.1.1) and (1.1.3).

In Section 4, we prove the main structural result, underlying the algorithm.

In Section 5, we compute theta functions associated with integer points in a subspace.

In Section 6, we estimate the distance from a given point to a lattice containing ${\Bbb Z}^n$.

In Section 7, we present an algorithm for sampling from a discrete Gaussian distribution.

In Section 8, we discuss the smooth case.

\head 3. The main algorithm \endhead

A function $G: {\Bbb R}^n \longrightarrow {\Bbb R}_+$ is called {\it log-concave} if 
$$G(\alpha x + (1-\alpha) y) \ \geq \ G^{\alpha}(x) G^{1-\alpha}(y) \quad \text{for all} \quad x, y \in {\Bbb R}^n \quad \text{and all} \quad 0 \leq \alpha \leq 1.$$
Equivalently, $G=e^{\psi}$ where $\psi: {\Bbb R}^n \longrightarrow {\Bbb R} \cup \{-\infty\}$ is concave, that is 
$$\psi\bigl(\alpha x + (1-\alpha) y\bigr) \ \geq \ \alpha \psi(x) + (1-\alpha) \psi(y) \quad \text{for all} \quad x, y \in {\Bbb R}^n \quad \text{and all} \quad 0 \leq \alpha \leq 1.$$

Recall that by $\|A\|_{\text{op}}$ we denote the operator norm of a matrix $A$, that is the largest singular value of $A$.

Our main result is as follows. 

\proclaim{(3.1) Theorem} Let $A=\left(a_{ij}\right)$ be an $m \times n$ real matrix, let $b=\left(\beta_1, \ldots, \beta_n\right)$ be a real $n$-vector and let $s >0$ be a real number. Let 
$$B=s I + {1 \over 2} A^T A$$
be an $n \times n$ positive definite matrix. 

Let $q=e^{-s}$ and
let us define a function 
$F_{A, b, s}: {\Bbb R}^m \longrightarrow {\Bbb R}_+$ by 
$$\split &F_{A, b, s}(t) = \prod_{j=1}^n \prod_{k=1}^{\infty} \left(1 + 2 q^{2k-1} \cos \left(\beta_j + \sum_{i=1}^m a_{ij} \tau_i \right) + q^{4k-2}\right), \\ &\qquad \text{where} \quad t=\left(\tau_1, \ldots, \tau_m\right). \endsplit$$
Then
\roster
\item We have 
$$\split &(2\pi)^{-m/2} \prod_{k=1}^{\infty} \left(1-q^{2k}\right)^n  \int_{{\Bbb R}^m} F_{A, b, s}(t) e^{-\|t\|^2/2} \ dt \\ 
&\qquad=\sum_{x \in {\Bbb Z}^n} \exp\left\{ - \langle B x, x \rangle +\ii \langle b, x \rangle \right\}. \endsplit$$
\item Suppose that 
$$\|A^T A\|_{\text{op}} \sum_{k=1}^{\infty} {q^{2k-1} \over \left(1-q^{2k-1}\right)^2} \ \leq \ {1 \over 2}.$$
Then for every  integer $K >0$ the function $G(t)=G_{A,b,s, K}(t)$ defined by
$$\split &G(t)= e^{-\|t\|^2/2} \prod_{j=1}^n \prod_{k=1}^K \left(1 + 2 q^{2k-1} \cos \left(\beta_j + \sum_{i=1}^m a_{ij} \tau_i \right) + q^{4k-2}\right),\\
&\qquad \text{where} \quad  t=\left(\tau_1, \ldots, \tau_m\right), \endsplit$$
is log-concave. In particular, the function $F_{A, b, s}(t) e^{-\|t\|^2/2}$ is log-concave. 
\endroster 
\endproclaim

We note that 
$$\sum_{k=1}^{\infty} {q^{2k-1} \over (1-q^{2k-1})^2} \ \leq \ {1 \over (1-q)^2} \sum_{k=1}^{\infty} q^{2k-1}={q \over (1-q)^2(1-q^2)} = {e^{-s} \over (1-e^{-s})^2(1-e^{-2s})}.$$
Consequently, to satisfy the constraint in Part (2), we are allowed to choose $A$ so that 
$$\|A^T A\|_{\text{op}} \ \leq \ {1 \over 2} e^s \left(1-e^{-s}\right)^2 \left(1-e^{-2s}\right).$$
We prove Theorem 3.1 in Section 4. 

Theorem 3.1 allows us to approximate $\Theta(B)$ and, more generally the sum (1.1.3), by using any of the efficient algorithms for integrating log-concave functions \cite{AK91}, \cite{F+94}, \cite{FK99},  \cite{LV07}. 

\subhead (3.2) Algorithm for computing the theta function \endsubhead 
We present an algorithm for computing (1.1.3).
\bigskip
{\bf Input:} An $n \times n$ positive definite matrix $B$ such that 
$$s I \ \preceq \ B \ \preceq \ \left(s + {e^s \over 4} \left(1-e^{-s}\right)^2\left(1-e^{-2s}\right)\right)I \quad \text{for some} \quad s \geq 1,$$
a vector $b \in {\Bbb R}^n$, $b=\left(\beta_1, \ldots, \beta_n\right)$, and a number $0 < \epsilon < 1$. 
\medskip
{\bf Output:} A positive real number approximating 
$$\sum_{x \in {\Bbb Z}^n} \exp\left\{ - \langle Bx, x \rangle +\ii \langle b, x \rangle \right\}$$
within relative error $\epsilon$.
\medskip
{\bf Algorithm:} Let $C=B-sI$. Hence $C$ is a positive definite matrix with 
$$\|C\|_{\text{op}} \ \leq \ {e^s \over 4} \left(1-e^{-s}\right)^2\left(1-e^{-2s}\right).$$
Next, we write 
$$C={1 \over 2} A^T A \quad \text{so that} \quad B=s I + {1 \over 2} A^T A$$
for an $m \times n$ matrix $A$. We can always choose $m=n$ or $m=\rk A$. Hence 
$$\left\| A^T A \right\|_{\text{op}} \ \leq \ {1 \over 2} e^s \left(1 - e^{-s}\right)^2 \left(1-e^{-2s}\right).$$ 
Let $q=e^{-s}$.
For an integer $K=K(\epsilon) >0$, to be specified in a moment, we define $\widehat{F}: {\Bbb R}^m \longrightarrow {\Bbb R}$ by 
$$\split &\widehat{F}(t)=\prod_{j=1}^n \prod_{k=1}^K \left(1 + 2 q^{2k-1} \cos \left(\beta_j + \sum_{i=1}^m a_{ij} \tau_i \right) + q^{4k-2}\right) 
\\ &\qquad \text{for} \quad t=\left(\tau_1, \ldots, \tau_m\right) \endsplit$$
and use any of the efficient algorithms of integration log-concave functions to compute 
$$(2\pi)^{-m/2} \prod_{k=1}^K \left(1-q^{2k}\right)^n \int_{{\Bbb R}^m} \widehat{F}(t) e^{-\|t\|^2/2} \ d t$$
within relative error $\epsilon/3$.

We choose $K$ so that the relative error acquired by replacing infinite products
$$\prod_{k=1}^{\infty} \left(1-q^{2k}\right)^n \quad \text{and} \quad \prod_{k=1}^{\infty} \left(1 + 2q^{2k-1} \cos\left(\beta_j +\sum_{i=1}^m a_{ij} \tau_i \right) 
+q^{4k-2}\right)$$
 in Theorem 3.1 by finite ones does not exceed $\epsilon/3$.
Since 
$$| \ln (1+x)| \ \leq \ 2|x| \quad \text{for} \quad -0.5 \ \leq \ x \ \leq \ 0.5,$$
and $q=e^{-s} \leq e^{-1}$, we have 
$$\left| \sum_{k=K}^{\infty} \ln \left(1-q^k\right)\right| \ \leq \ 2 \sum_{k=K}^{\infty} q^k = {2 q^K \over 1-q} \ \leq 4 q^K.$$
Similarly, 
$$\split &\left| \sum_{k=K}^{\infty} \ln \left(1+ 2q^{2k-1} \cos\left(\beta_j + \sum_{i=1}^m a_{ij} \tau_i \right) + q^{4k-2}\right) \right| \\ &\qquad \leq \
\left| \sum_{k=K}^{\infty} \ln \left(1 -2 q^{2k-1} + q^{4k-2} \right) \right|  = 2 \left| \sum_{k=K}^{\infty} \ln \left(1-q^{2k-1}\right) \right| \\
&\qquad \leq 4 \sum_{k=K}^{\infty}  q^{2k-1} ={4 q^{2K-1} \over 1-q^2} \ \leq \ 5 q^{2K-1}. \endsplit$$
Consequently, to approximate the infinite products in Theorem 3.1 by finite ones within relative error $\epsilon/3$, we can choose 
$K=O\left(\ln (n/\epsilon)\right)$. 
We summarize the result as a theorem.
\proclaim{(3.3) Theorem} Given an $n \times n$ positive definite matrix $B$ satisfying (2.1.1), a vector $b \in {\Bbb R}^n$ and $0 < \epsilon \leq 1$, the algorithm of Section 3.2 approximates 
$$\sum_{x \in {\Bbb Z}^n} \exp\left\{ -\langle Bx, x \rangle + \ii \langle b, x \rangle \right\}$$
within relative error $\epsilon$ 
in time polynomial in $n$, $s$ and $\epsilon^{-1}$.
\endproclaim
{\hfill \hfill \hfill} \qed

\head 4. Proof of Theorem 3.1 \endhead 

The proof of Part (1) is based on the Jacobi identity.
\subhead (4.1) Jacobi's formula \endsubhead  For any $0 \leq q < 1$ and any $w \in {\Bbb C} \setminus 0$, we have 
$$\prod_{k \geq 1} \left(1-q^{2k}\right) \left(1+ w q^{2k-1}\right) \left(1+ w^{-1} q^{2k-1}\right) = \sum_{\xi \in {\Bbb Z}} w^{\xi} q^{\xi^2}. $$
This is Jacobi's triple product identity, see for example, Section 2.2 of \cite{An98}. Suppose now that 
$$w_j \in {\Bbb C} \setminus \{0\} \quad \text{for} \quad j=1, \ldots, n.$$
Then 
$$\aligned & \prod_{j=1}^n \prod_{k \geq 1} \left(1-q^{2k}\right) \left(1+ w_j q^{2k-1}\right) \left(1+ w_j^{-1} q^{2k-1}\right) \\ &\qquad =
\sum\Sb x\in {\Bbb Z}^n: \\ x =\left(\xi_1, \ldots, \xi_n\right) \endSb q^{\|x\|^2} \prod_{j=1}^n w_j^{\xi_j}. \endaligned \tag4.1.1$$

\subhead (4.2) Proof of Part (1) \endsubhead For $t=\left( \tau_1, \ldots, \tau_m \right)$, we choose 
$$w_j(t) =\exp\left\{ \ii \left(\beta_j + \sum_{i=1}^m a_{ij} \tau_i \right)\right\} \quad \text{for} \quad j=1, \ldots, n$$
in (4.1.1). Using that 
$$\split &\left(1+ w_j(t) q^{2k-1}\right)\left(1+ w_j^{-1}(t) q^{2k-1}\right)=1+ \left(w_j(t) + w_j^{-1}(t) \right) q^{2k-1} + q^{4k-2} \\&\qquad=
1 + 2 \cos\left(\beta_j + \sum_{i=1}^m a_{ij} \tau_i \right) q^{2k-1} + q^{4k-2} \endsplit$$
and that 
$$\prod_{j=1}^n w_j^{\xi_j} =\exp\left\{ \ii \sum_{j=1}^n \beta_j \xi_j + \ii \sum_{i=1}^m \tau_i  \left(\sum_{j=1}^n a_{ij} \xi_j \right)\right\},$$
we conclude that 
$$\split &F_{A, b, s}(t) \prod_{k=1}^{\infty} \left(1-q^{2k}\right)^n \\&\qquad = \sum\Sb x \in {\Bbb Z}^n: \\ x=\left(\xi_1, \ldots, \xi_n\right) \endSb 
q^{\|x\|^2} \exp\left\{ \ii \sum_{j=1}^n \beta_j \xi_j + \ii \sum_{i=1}^m \tau_i  \left( \sum_{j=1}^n a_{ij} \xi_j \right) \right\}. \endsplit$$
Since 
$${1 \over \sqrt{2 \pi}} \int_{-\infty}^{+\infty} \exp\left\{ \ii \tau_i \sum_{j=1}^n a_{ij} \xi_j \right\} e^{-\tau_i^2/2} \ d \tau_i = 
\exp\left\{ -{1 \over 2} \left(\sum_{j=1}^n a_{ij} \xi_j \right)^2 \right\},$$
we get 
$$\split &(2 \pi)^{-m/2} \prod_{k=1}^{\infty} \left(1-q^{2k}\right)^n \int_{{\Bbb R}^m} F_{A, b, s}(t) e^{-\|t\|^2/2} \ dt \\&\qquad =
\sum\Sb x \in {\Bbb Z}^n: \\ x =\left(\xi_1, \ldots, \xi_n\right)\endSb q^{\|x\|^2} \exp\left\{ - {1 \over 2} \sum_{i=1}^m \left( \sum_{j=1}^n a_{ij} \xi_j \right)^2
+\ii \sum_{j=1}^n \beta_j \xi_j  \right\} \\
&\qquad =\sum_{x \in {\Bbb Z}^n} q^{\|x\|^2} \exp\left\{ -{1 \over 2} \|A x\|^2 + \ii \langle b, x \rangle \right\} = \sum_{x \in {\Bbb Z}^n} \exp\left\{ - \langle Bx , x \rangle + \ii \langle b, x \rangle \right\},
\endsplit$$
and the proof follows.
{\hfill \hfill \hfill} \qed

To prove Part (2), we need one technical estimate.

\proclaim{(4.3) Lemma} Let $0 < q < 1$ and $\alpha, \beta$ be reals. Then 
$${d^2 \over d\tau^2} \ln \left( 1+ 2q\cos(\alpha \tau + \beta) + q^2\right) \ \leq \ {2 \alpha^2 q \over (1-q)^2}.$$
\endproclaim 
\demo{Proof}
We have 
$${d \over d\tau} \ln \left(1 + 2q \cos (\alpha \tau + \beta)  + q^2\right)=
-{2 \alpha q \sin (\alpha \tau + \beta) \over 1 + 2 q\cos (\alpha \tau+ \beta)  + q^2}$$
and 
$$\split &{d^2 \over d \tau^2} \ln \left(1 + 2 q\cos (\alpha \tau  + \beta)  + q^2\right)\\=&\qquad -{2\alpha^2 q \cos(\alpha \tau +\beta) \left(1+ 2 q \cos(\alpha \tau+\beta)  +q^2\right) +
\left(2 \alpha q \sin(\alpha \tau +\beta) \right)^2
 \over \left(1 +2 q\cos (\alpha \tau + \beta) +q^2 \right)^2} \\
 =&-{ 2 \alpha^2 q (1+q^2) \cos (\alpha \tau + \beta) + 4\alpha^2 q^2 \over  \left(1 +2 q\cos (\alpha \tau + \beta) +q^2 \right)^2}.  \endsplit$$
 Now, 
 $$\left(1 + 2q \cos(\alpha \tau+\beta) + q^2\right)^2 \ \geq \ \left(1-2q+q^2\right)^2=(1-q)^4.$$
 Also,
 $$\split &2 \alpha^2 q (1+q^2)\cos (\alpha \tau +\beta)  + 4 \alpha^2 q^2 \ \geq \  -2\alpha^2 q (1+q^2) + 4\alpha^2 q^2 \\&\qquad = 2\alpha^2 q\left(2 q -1-q^2\right)  =  -2 \alpha^2 q (1-q)^2. \endsplit$$
The proof now follows.
 {\hfill \hfill \hfill} \qed
 \enddemo
 
 \subhead (4.4) Proof of Part (2) \endsubhead
 It suffices to prove that the restriction of $G(t)$ onto any affine line 
 $$\tau_i = \gamma_i \tau + \delta_i \quad \text{for} \quad i=1, \ldots, m \quad \text{where} \quad \sum_{i=1}^m \gamma_i^2 =1$$
 is log-concave. Indeed, let $g(\tau)$ be that restriction. From Lemma 4.3, we get 
 $$\split {d^2 \over d \tau^2} \ln g(\tau) \ \leq \ &-1 + 2 \sum_{k=1}^K  {q^{2k-1} \over \left(1-q^{2k-1}\right)^2} \sum_{j=1}^n \left(\sum_{i=1}^m a_{ij} \gamma_i \right)^2 \\
 &\leq -1 + 2 \|A^T\|^2_{\text{op}} \sum_{k=1}^K {q^{2k-1} \over \left(1-q^{2k-1}\right)^2} \\
 &=-1 + 2 \|A^T A \|_{\text{op}} \sum_{k=1}^K {q^{2k-1} \over \left(1 -q^{2k-1}\right)^2} \ \leq \ 0 \endsplit$$
 and hence $\ln g(\tau)$ is concave. The proof now follows.
 {\hfill \hfill \hfill} \qed 
 
 \head 5. Integer points in a subspace \endhead
 
 Let $A$ be an $m \times n$ integer matrix of $\rk A = m < n$ and let $L=\ker A$ be a subspace, $L \subset {\Bbb R}^n$. Then 
 $\Lambda ={\Bbb Z}^n \cap L$ is a lattice in $L$. Note that in this case, we do not define $\Lambda$ by its basis.
 For $s >0$, we consider the theta function
 $$\Theta_{\Lambda}(s)=\sum_{x \in \Lambda} e^{-s \|x\|^2}.$$
 Our main result is as follows.
 
 \proclaim{(5.1) Theorem} 
Suppose that $\|A\|_{\text{op}} \leq \gamma$ for some $\gamma \geq 1$. 
For $s > 0$ and $t >0$, let $B=B_{s,t}$ be an $n \times n$ positive definite matrix with the eigenvectors in $L \cup L^{\bot}$, where $L=\ker A$,  and such that the eigenvectors in 
$L$ have eigenvalue $s$ while the eigenvectors in $L^{\bot}$ have eigenvalue $s+t$.  Then 
$$\left| \Theta(B)-\Theta_{\Lambda}(s) \right| \ \leq \ \exp\left\{ -{t \over  \gamma^2} + {2n e^{-s} \over 1-e^{-s}} \right\}.$$
\endproclaim 

\example{(5.2) Example}
Let us fix $\delta >0$ and let 
$$s=\left({1 \over 2} + \delta \right) \ln n \quad \text{and} \quad t={e^s \over 5}={n^{{1 \over 2}+\delta} \over 5}.  \tag5.2.1$$ 
From Theorem 5.1, we have 
$$\left| \Theta(B) - \Theta_{\Lambda}(s)\right| \ \leq \ \exp\left\{ -{n^{{1 \over 2} + \delta} \over 5 \gamma^2} + {2n^{{1\over 2} -\delta} \over 
1-n^{-{1 \over 2}-\delta}} \right\}.$$
As long as $\gamma=o\left(n^{\delta}\right)$, we get 
$$\left| \Theta(B)-\Theta_{\Lambda}(s) \right| =o(1). \tag5.2.2$$
When $s \geq 3$, the matrix $B=B_{s,t}$ satisfies (2.1.2) and hence $\Theta(B)$ can be efficiently approximated. Since $\Theta(B) \geq 1$, from (5.2.2) and Theorem 3.3, we obtain a randomized polynomial time algorithm that approximates $\Theta_{\Lambda}(s)$ within a relative error of $o(1)$ as $n \longrightarrow \infty$.
\endexample

The proof of Theorem 5.1 is based on the following two lemmas. In the first lemma, we bound from below the distance of a point $x \in {\Bbb Z}^n \setminus 
\Lambda$ to the subspace $L$.
\proclaim{(5.3) Lemma} Let $A$ be an $m \times n$ integer matrix with $\rk A = m < n$ and let $L=\ker A$. For a point 
$x \in {\Bbb R}^n$, let 
$$\dist(x, L)=\min_{y \in L} \|x-y\|$$ be the Euclidean distance from $x$ to $L$. 
Then 
$$\dist(x, L) \ \geq \ \left( \|A\|_{\text{op}} \right)^{-1} \quad \text{for all} \quad x \in {\Bbb Z}^n \setminus L.$$
\endproclaim
\demo{Proof} Suppose that $x \in {\Bbb Z}^n \setminus L$. Let $P: {\Bbb R}^n \longrightarrow L^{\bot} = \operatorname{image} A^T$ be the orthogonal projection. Then the matrix of $P$ in the standard coordinates is $A^T (A A^T)^{-1} A$ and hence
$$\dist^2 (x, L) =\left\| P(x)\right\|^2 = \big\langle A^T (A A^T)^{-1} Ax, \ A^T (A A^T)^{-1} Ax \big\rangle =
\big\langle (A A^T)^{-1} Ax, \ Ax \big\rangle.$$
Since $A$ is an integer matrix, $x$ is an integer vector and $Ax \ne 0$, we have $\|Ax\| \geq 1$. Let $\lambda > 0$ be the smallest eigenvalue 
of the matrix $(A A^T)^{-1}$. Then 
$$\big\langle (A A^T)^{-1} Ax, \ Ax \big\rangle \ \geq \ \lambda \|Ax\|^2 \ \geq \ \lambda$$
and hence 
$$\dist^2(x, L) \ \geq \ \lambda.$$
On the other hand,
$$\lambda = \left( \| A A^T \|_{\text{op}}\right)^{-1} =\left( \|A\|_{\text{op}}\right)^{-2},$$
from which the proof follows.
{\hfill \hfill \hfill} \qed
\enddemo 

The next lemma provides some technical estimates for the theta function. For the proof of Theorem 5.1 we need Part (1) only, while Part (2) will be used later.
\proclaim{(5.4) Lemma} 
\roster
\item For $s >0$, we have 
$$\Theta(sI) =\sum_{x \in {\Bbb Z}^n} e^{-s\|x\|^2} \ \leq \ \exp\left\{ {2n e^{-s} \over 1-e^{-s}}\right\}.$$
\item For $s>0$ and 
$$4n e^{-1} \ \geq \ k \ \geq \ 30 n e^{-s},$$
we have 
$$\sum\Sb x \in {\Bbb Z}^n: \\ \|x\|^2 \geq k \endSb e^{-s\|x\|^2} \ \leq \ e^{-k}.$$
\endroster
\endproclaim
\demo{Proof}
For $s > 0$, we have 
$$\split \Theta(sI)=  &\left(\sum_{\xi \in {\Bbb Z}} e^{-s\xi^2} \right)^n  \leq \ \left(1 + 2\sum_{\xi=1}^{\infty} e^{-s \xi} \right)^n =
\left(1 + {2 e^{-s} \over 1-e^{-s}}\right)^n \\=&\exp\left\{ n \ln \left(1 + {2 e^{-s} \over 1-e^{-s}}\right)\right\} \ \leq \ \exp\left\{{2n e^{-s} \over 1-e^{-s}}\right\},
\endsplit$$
which proves Part (1).

To prove Part (2), for any $0 < \tau < s$, using Part (1), we get 
$$\split \sum\Sb x \in {\Bbb Z}^n: \\ \|x\|^2 \geq k \endSb e^{-s\|x\|^2} \ \leq \ &e^{-\tau k} \sum\Sb x \in {\Bbb Z}^n: \\ \|x\|^2 \geq k \endSb e^{-s\|x\|^2} e^{\tau \|x\|^2} 
\ \leq \ e^{-\tau k} \Theta\bigl((s-\tau) I \bigr) \\ \leq \ &\exp\left\{-\tau k + {2n e^{-(s-\tau)} \over 1-e^{-(s-\tau)}} \right\}. \endsplit$$
Optimizing on $\tau$, we choose 
$$\tau=s + \ln {k \over 4n}.$$
Since $k \geq 30 ne^{-s}$, we have
$$\tau \ \geq \ \ln {30 \over 4} \ > \ 2$$
and since $k \leq 4n e^{-1}$, we have 
$$s-\tau = -\ln {k \over 4n} \ \geq \ 1.$$
Therefore, 
$$\sum\Sb x \in {\Bbb Z}^n: \\ \|x\|^2 \geq k \endSb e^{-s\|x\|^2} \ \leq \ \exp\left\{- \tau k + 4n e^{-(s-\tau)}\right\} = \exp\left\{ -(\tau-1)k\right\} \ \leq \ e^{-k},$$
as required.
{\hfill \hfill \hfill} \qed
\enddemo

Now we are ready to prove Theorem 5.1.

\subhead (5.5) Proof of Theorem 5.1 \endsubhead 

Applying Lemma 5.3 and  Part(1) of  Lemma 5.4, we obtain
$$\split \left| \Theta(B) - \Theta_{\Lambda}(s)\right|  = &\sum_{x \in {\Bbb Z^n} \setminus L} \exp\left\{ - \langle Bx, x \rangle \right\} \\ = \ 
&\sum_{x \in {\Bbb Z}^n \setminus L} \exp\left\{-t \dist^2(x, L)\right\} \exp\left\{-s\|x\|^2\right\} \\ 
\leq \ &\exp\left\{ -{t \over  \gamma^2}\right\} \sum_{x \in {\Bbb Z}^n} \exp\left\{-s \|x\|^2\right\} \\ 
\leq \ &\exp\left\{ - {t \over \gamma^2}  + {2n e^{-s} \over 1-e^{-s}}\right\}.
\endsplit$$
{\hfill \hfill \hfill} \qed

As in Example 5.2, let us fix $0 < \delta < {1 \over 2}$, define $s$ and $t$ by (5.2.1) and assume that $\|A\|_{\text{op}} =o\left(n^{\delta}\right)$, so that 
$\Theta_{\Lambda}(s)$ can be approximated in randomized polynomial time within a relative error of $o(1)$.
 If there are no points $x \in \Lambda \setminus \{0\}$ with $\|x\|^2 \leq 30n^{{1 \over 2}-\delta}$ then by Part (2) of Lemma 5.4, we have $\Theta_{\Lambda}(s) =1+o(1)$. On the other hand, 
 if $\Lambda$ contains many short vectors, then $\Theta_{\Lambda}(s)$ can be large. For example, if $L$ is a coordinate subspace, $\dim L \geq \alpha n$ for 
 some $0 < \alpha < 1$, 
 so that $\Lambda$ is identified with ${\Bbb Z}^{\dim L}$, then 
 $$\Theta_{\Lambda}(s) \ \geq \ \left(\sum_{\xi \in {\Bbb Z}} e^{-s\xi^2}\right)^{\alpha n} \ \geq \ \left(1 + 2e^{-s}\right)^{\alpha n}=\left(1 + {2 \over n^{{1 \over 2}+\delta}}\right)^{\alpha n} \ \geq \ \exp\left\{ \alpha n^{{1 \over 2}-\delta}\right\}$$
 is exponentially large in $n$. Hence computing $\Theta_{\Lambda}(s)$ allows us to distinguish the case of $L$ having no short non-zero integer vectors from the case of $L$ having sufficiently many short integer vectors.
 
\head 6. Lattices containing ${\Bbb Z}^n$ \endhead

As in Section 5, for a lattice $\Lambda \subset {\Bbb R}^n$, a point $v \in {\Bbb R}^n$ and a number $\tau > 0$, we denote 
$$\Theta_{\Lambda}(\tau, v)=\sum_{u \in \Lambda} \exp\left\{ -\tau \|u-v\|^2 \right\}.$$
In agreement with our notation in Sections 1-4, when $\Lambda = {\Bbb Z}^n$, we still denote 
$\Theta_{{\Bbb Z}^n}(\tau, v)$ just by $\Theta(\tau I, v)$ and $\Theta_{{\Bbb Z}^n}(\tau, 0)$ just by $\Theta(\tau I)$, so
$$\Theta(\tau I, v) =\sum_{x \in {\Bbb Z}^n} e^{-\tau \|x-v\|^2} \quad \text{and} \quad \Theta(\tau I)=\sum_{x \in {\Bbb Z}^n} e^{-\tau \|x\|^2}.$$
In this section we prove the following main result. 
\proclaim{(6.1) Theorem} Let $\Lambda \subset {\Bbb R}^n$ be a lattice such that ${\Bbb Z}^n \subset \Lambda$. Then for $0 < \tau \leq  1$, we have 
$$41 e^{-\pi^2/\tau} \dist^2(v, \Lambda) \ \geq \ \ln {\Theta(\tau I) \over \Theta_{\Lambda}(\tau, v)} \ \geq \ 
13 e^{-\pi^2/\tau} \dist^2(v, \Lambda) + \ln \det \Lambda.$$
\endproclaim

Apart from the additive term of $\ln \det \Lambda$, the formula of Theorem 6.1 provides an estimate of 
$\dist(v, \Lambda)$ within a constant factor of $\sqrt{41/13} \approx 1.8$. It may happen that as $n$ grows, the additive term becomes asymptotically negligible, and hence the formula of Theorem 6.1 provides an approximation of $\dist(v, \Lambda)$ within a constant factor. 

\example{(6.2) Example}
A lattice $\Lambda \subset {\Bbb R}^n$ containing ${\Bbb Z}^n$ can be constructed as follows: let $w_1, \ldots, w_n$ be a basis of ${\Bbb Z}^n$ and let 
$\lambda_1, \ldots, \lambda_n$ be positive integers. Then 
$$u_i= {1 \over \lambda_i} w_i \quad \text{for} \quad i=1, \ldots, n \tag6.2.1$$
is a basis of a lattice $\Lambda$ containing ${\Bbb Z}^n$. Moreover, any lattice containing ${\Bbb Z}^n$ can be constructed this way, cf., for example, Chapter I of \cite{Ca97} for the Smith normal form. We have 
$$\ln \det \Lambda = - \sum_{i=1}^n \ln \lambda_i.$$
Let us consider the case when 
$\dist^2(v, \Lambda) \geq n^{\alpha}$ for some $0 < \alpha < 1$.
 We let
 $$\tau={10\pi^2 \over \alpha \ln n},$$ so that 
$$e^{-\pi^2/\tau} \dist^2(v, \Lambda)= n^{-0.1 \alpha} \dist^2(v, \Lambda)\ \geq \ n^{0.9 \alpha}.$$
To make sure that the term $\ln \det \Lambda$ is asymptotically negligible, we choose not more than $n^{0.8 \alpha}$ of $\lambda_i$ in (6.2.1) satisfying 
$\lambda_i \leq \gamma$ for a constant $\gamma > 1$, fixed in advance, while the rest of $\lambda_i$ are equal to 1. 

Let $B$ be the Gram matrix of the basis $u_1, \ldots, u_n$. In the trivial case, if $w_1, \ldots, w_n$ in (6.2.1) is the standard basis $e_1, \ldots, e_n$, then for large $n$, the matrix $\tau B$ satisfies (2.1.4) and hence the ratio $\Theta(\tau I)/\Theta_{\Lambda}(\tau, v) $ can be approximated in randomized polynomial time. However, the matrix $\tau B$ would still satisfy (2.1.4) in a less trivial situation, when $w_1, \ldots, w_n$ are close enough to the standard basis, for example, when $w_i =A e_i$ for some matrix $A \in GL(n, {\Bbb Z})$ where 
$$ \|A\|_{\text{op}} \ \leq \ \gamma \quad \text{and} \quad \|A^{-1}\|_{\text{op}} \ \leq \ n^{\alpha/21 \gamma^2},$$
for a constant $\gamma > 1$, fixed in advance.

It appears essential that we are able to choose $\tau$ in the non-smooth range, see Section 2.1.5. Indeed, choosing $\tau \leq \pi^2/\gamma \ln n$ for some 
$\gamma > 1$ leads to 
$$e^{-\pi^2/\tau} \dist^2(v, \Lambda)=o(1)$$ and 
hence the $\ln \det \Lambda$ additive term cannot be discarded. 
\endexample
We note that the ratio $\Theta_{\Lambda}(\tau, v)/\Theta_{\Lambda}(\tau, 0)$ was crucially used by Aharonov and Regev to show that estimating $\dist(v, \Lambda)$ within a factor of $O(\sqrt{n})$ for any lattice $\Lambda \subset {\Bbb R}^n$ lies in NP $\cap$ co-NP \cite{AR05}.

To prove Theorem 6.1, we first consider the case of $\Lambda={\Bbb Z}^n$.

 \proclaim{(6.3) Lemma} For $y \in {\Bbb R}^n$ and $0 < \tau \leq 1$, we have 
 $$\exp\left\{ -41 e^{-\pi^2/\tau} \dist^2(y, {\Bbb Z}^n)\right\} \leq \ { \Theta(\tau I, y)\over \Theta(\tau I)} \ \leq \ \exp\left\{ -13 e^{-\pi^2/\tau} \dist^2\left(y, {\Bbb Z}^n\right)\right\}.$$
\endproclaim
\demo{Proof} Let $y=\left(\eta_1, \ldots, \eta_n \right)$. We have 
$$\Theta(\tau I, y) =\sum_{x \in {\Bbb Z}^n} \exp\left\{ - \tau\|x-y\|^2 \right\} = \prod_{i=1}^n \sum_{\xi \in {\Bbb Z}} \exp\left\{ - \tau (\xi-\eta_i)^2 \right\} $$
and similarly,
$$\Theta(\tau I)=\sum_{x \in {\Bbb Z}^n} \exp\left\{ - \tau \|x\|^2 \right\} = \prod_{i=1}^n \sum_{\xi \in {\Bbb Z}} \exp\left\{- \tau \xi^2 \right\}.$$
Translating $y$ by an integer vector, without loss of generality we assume that $y=\left(\eta_1, \ldots, \eta_n \right)$ where 
$$|\eta_i| \ \leq \ {1 \over 2} \quad \text{for} \quad i=1, \ldots, n.$$
Then
$$\dist^2(y, {\Bbb Z}^n) = \|y\|^2=\sum_{i=1}^n \eta_i^2.$$

By the reciprocity relation (1.1.4), we get 
$$\split &\Theta(\tau I, y) ={\pi^{n/2} \over \tau^{n/2}} \prod_{i=1}^n \sum_{\xi \in {\Bbb Z}} \exp\left\{ - \pi^2 \tau^{-1} \xi^2 + 2 \pi \ii \xi \eta_i\right\} \quad \text{and} \\
&\Theta(\tau I) ={\pi^{n/2} \over \tau^{n/2}} \prod_{i=1}^n \sum_{\xi \in {\Bbb Z}} \exp\left\{ -\pi^2 \tau^{-1} \xi^2 \right\}. \endsplit$$
Denoting 
$$q=e^{-\pi^2/\tau},$$
from the Jacobi identity (4.1), we get 
$$\split &\sum_{\xi \in {\Bbb Z}} \exp\left\{-\pi^2 \tau^{-1} \xi^2 + 2 \pi \ii \xi \eta_i\right\} \\&\qquad = 
\prod_{k=1}^{\infty} \left(1-q^{2k}\right)\left(1 + \exp\left\{2 \pi \ii  \eta_i \right\} q^{2k-1}\right) \left(1+ \exp\left\{ -2 \pi \ii  \eta_i \right\} q^{2k-1}\right)
\\ &\qquad =\prod_{k=1}^{\infty} \left(1-q^{2k}\right) \left(1 + 2q^{2k-1}  \cos (2 \pi \eta_i)  + q^{4k-2}\right) \endsplit$$
and, similarly,
$$\sum_{\xi \in {\Bbb Z}} \exp\left\{-\pi^2 \tau^{-1} \xi^2 \right\}=\prod_{k=1}^{\infty}\left(1-q^{2k}\right) \left(1 + 2 q^{2k-1} + q^{4k-2}\right).$$
Summarizing,
$$\split {\Theta(\tau I, y) \over \Theta(\tau I)} =& \prod_{i=1}^n \prod_{k=1}^{\infty} {1+ 2 q^{2k-1} \cos( 2 \pi \eta_i) +q^{4k-2} \over 1+ 2q^{2k-1} +q^{4k-2}} \\
=&\prod_{i=1}^n \prod_{k=1}^{\infty} \left(1 - {2 q^{2k-1} \left(1-\cos(2 \pi \eta_i)\right) \over \left(1+q^{2k-1}\right)^2}\right). \endsplit$$
We have 
$$7 \eta^2 \ \leq \ 1-\cos(2 \pi \eta) \ \leq \ 20 \eta^2 \quad \text{for} \quad -{1 \over 2} \leq \eta \leq {1 \over 2}.$$
Since 
$$q=e^{-\pi^2/\tau} \ \leq \ e^{-\pi^2} \ < \ 10^{-4} \quad \text{and} \quad |\eta_i| \ \leq \ {1 \over 2},$$
we have 
$$ {\eta_i^2 q^{2k-1} \over \left(1+q^{2k-1}\right)^2} \ \leq \ {1 \over 4} 10^{-4}$$
and we can further write 
$$\prod_{i=1}^n \prod_{k=1}^{\infty} \left(1 - {40\eta_i^2 q^{2k-1} \over \left(1+q^{2k-1}\right)^2}\right) \ \leq \ {\Theta(\tau I, y) \over \Theta(\tau I)} \ \leq \ \prod_{i=1}^n \prod_{k=1}^{\infty} \left(1 - {14\eta_i^2 q^{2k-1} \over \left(1+q^{2k-1}\right)^2}\right) \tag6.3.1$$
(note that all factors in the products are positive).
 
 Using that  
$$\ln (1- \alpha) \ \leq \ -\alpha \quad \text{for} \quad 0 \leq \alpha < 1,$$
we conclude that 
$$\aligned &\prod_{k=1}^{\infty} \left(1-{14 \eta_i^2 q^{2k-1} \over (1+q^{2k-1})^2}\right) = 
\exp\left\{ \sum_{k=1}^{\infty} \ln\left(1- {14 \eta_i^2 q^{2k-1} \over (1+q^{2k-1})^2}\right)\right\} \\ &\qquad \leq \ \exp\left\{ -\sum_{k=1}^{\infty} 
{14 \eta_i^2 q^{2k-1} \over (1+q^{2k-1})^2}\right\} \ \leq \ \exp\left\{ -13 \eta_i^2 \sum_{k=1}^{\infty} q^{2k-1}\right\} 
\\ &\qquad = \exp\left\{ - {13 \eta_i^2 q \over 1-q^2}\right\} \ \leq \ \exp\left\{ -13 \eta_i^2 q \right\}. \endaligned \tag6.3.2$$
Similarly, using that 
$$\ln (1-\alpha) \ \geq \ -1.01 \alpha \quad \text{for} \quad 0 \leq \alpha \leq 0.001,$$
we conclude that 
$$\aligned &\prod_{k=1}^{\infty} \left(1-{40 \eta_i^2 q^{2k-1} \over (1+q^{2k-1})^2}\right) = 
\exp\left\{ \sum_{k=1}^{\infty} \ln\left(1- {40 \eta_i^2 q^{2k-1} \over (1+q^{2k-1})^2}\right)\right\} \\ &\qquad \geq \ \exp\left\{ -\sum_{k=1}^{\infty} 
{40.5 \eta_i^2 q^{2k-1} \over (1+q^{2k-1})^2}\right\} \ \geq \ \exp\left\{ -40.5 \eta_i^2 \sum_{k=1}^{\infty} q^{2k-1}\right\} 
\\ &\qquad = \exp\left\{ - {40.5 \eta_i^2 q \over 1-q^2}\right\} \ \geq \ \exp\left\{ -41 \eta_i^2 q \right\}. \endaligned \tag6.3.3$$
Summarizing, from (6.3.1)--(6.3.3) we infer that 
$${\Theta(\tau I, y) \over \Theta(\tau I)} \ \leq \ \prod_{i=1}^n \exp\left\{ -13 \eta_i^2 q\right\} = \exp\left\{ -13 q \sum_{i=1}^n \eta_i^2\right\} 
= \exp\left\{-13 q \dist^2\left(y, {\Bbb Z}^n\right)\right\}$$
and 
$${\Theta(\tau I, y) \over \Theta(\tau I)} \ \geq \ \prod_{i=1}^n \exp\left\{ -41 \eta_i^2 q\right\} = \exp\left\{ -41 q \sum_{i=1}^n \eta_i^2\right\} 
= \exp\left\{-41 q \dist^2\left(y, {\Bbb Z}^n\right)\right\},$$
which concludes the proof.
{\hfill \hfill \hfill} \qed
\enddemo

Now we are ready to prove Theorem 6.1.
\subhead (6.4) Proof of Theorem 6.1 \endsubhead 
 Let $u_i, i \in I$ be the coset representatives of the quotient $\Lambda/{\Bbb Z}^n$, so that $\Lambda$ is represented as a disjoint union
$$\Lambda= \bigcup_{i \in I} \left(u_i + {\Bbb Z}^n\right) \quad \text{and} \quad |I| = {1 \over \det \Lambda}. \tag6.4.1$$
Then
$$\aligned \Theta_{\Lambda}(\tau, v)= &\sum_{u \in \Lambda} \exp\left\{ - \tau \|u-v\|^2\right\} = \sum_{i \in I} \sum_{x \in {\Bbb Z}^n} 
\exp\left\{ -\tau \|u_i + x-v \|^2\right\}\\=& \sum_{i \in I} \Theta(\tau I, v-u_i). \endaligned \tag6.4.2$$
On the other hand, 
$$\dist(v, \Lambda) = \min_{i \in I} \dist\left(v, u_i + {\Bbb Z}^n\right)=\min_{i \in I} \dist\left(v-u_i, {\Bbb Z}^n\right).$$
By Lemma 6.3, we have 
$$\Theta(\tau I, v-u_i) \ \leq \ \exp\left\{ -13 e^{-\pi^2/\tau} \dist^2\left(v-u_i, {\Bbb Z}^n\right)\right\} \Theta(\tau I)$$
and hence 
$$\Theta(\tau I, v-u_i) \ \leq \ \exp\left\{ -13 e^{-\pi^2/\tau} \dist^2(v, \Lambda) \right\} \Theta(\tau I). $$
Therefore, by (6.4.2)  we have 
$$\Theta_{\Lambda}(\tau, v) \ \leq \ |I| \exp\left\{ -13 e^{-\pi^2/\tau} \dist^2(v, \Lambda) \right\} \Theta(\tau I)$$
and from (6.4.1) we obtain
$${\Theta_{\Lambda}(\tau, v) \over \Theta(\tau I)} \ \leq \ (\det \Lambda)^{-1} \exp\left\{ -13 e^{-\pi^2/\tau} \dist^2(v, \Lambda) \right\}. \tag6.4.3$$
We have 
$$\dist(v, \Lambda) =\dist(v-u_{i_0}, {\Bbb Z}^n) \quad \text{for some} \quad i_0 \in I.$$
Therefore, by Lemma 6.3, 
$$\Theta(\tau I, v-u_{i_0}) \ \geq \ \exp\left\{ -41 e^{-\pi^2/\tau}  \dist^2(v, \Lambda)\right\} \Theta (\tau I).$$ 
Hence by (6.4.2)
$${\Theta_{\Lambda}(\tau, v) \over \Theta(\tau I)} \ \geq \ \exp\left\{ -41 e^{-\pi^2/\tau}  \dist^2(v, \Lambda)\right\}. \tag6.4.4$$
Combining (6.4.3)--(6.4.4), we complete the proof.
{\hfill \hfill \hfill \hfill} \qed

 \head 7. Sampling from the discrete Gaussian measure \endhead
 
 \subhead (7.1) Gaussian measure on lattices \endsubhead 
 Let $\Lambda \subset {\Bbb R}^n$ be a lattice and let $v \in {\Bbb R}^n$. In this section, we use the shorthand 
 $$\Theta_{\Lambda}(v)=\sum_{u \in \Lambda} e^{-\|u -v \|^2} \quad \text{and} \quad \Theta_{\Lambda}(0)=\sum_{u \in \Lambda} e^{-\|u\|^2}.$$
We consider the discrete Gaussian probability measure on $\Lambda$ defined by 
 $$\Pr(u)={\exp\{-\|u-v\|^2 \} \over \Theta_{\Lambda}(v)} \quad \text{for} \quad u \in \Lambda. \tag7.1.1$$
 Our goal is to sample a point $u \in \Lambda$ from a probability distribution that is $\epsilon$-close in the total variation distance to (7.1.1).
 
 Let $u_1, \ldots, u_n$ be a basis of $\Lambda$, so that every point $u \in {\Bbb R}^n$ can be uniquely written as
 $$u=\xi_1 u_1 + \ldots + \xi_n u_n \quad \text{for some} \quad \xi_1, \ldots, \xi_n \in {\Bbb R}, \tag7.1.2$$
 and $u \in \Lambda$ if and only if $\xi_1, \ldots, \xi_n$ are integer.
 
 The general design of the algorithm is the same as in \cite{G+08} and \cite{Pe10}: we consecutively sample the coordinates $\xi_n, \xi_{n-1}, 
 \ldots, \xi_1$ of $u$. For that, we compute the conditional distribution of $\xi_{n-k}$ for fixed $\xi_n, \ldots, \xi_{n-k+1}$.

 For $\alpha \in {\Bbb Z}$, let $H_{\alpha} \subset {\Bbb R}^n$ be the affine hyperplane defined by the equation 
 $\xi_n = \alpha$ in (7.1.2). Let $\Lambda_{\alpha} =\Lambda \cap H_{\alpha}$. We identify $H_{\alpha}$ with ${\Bbb R}^{n-1}$ by choosing the origin at a point of $\Lambda_{\alpha}$, so that $\Lambda_{\alpha} \subset H_{\alpha}$ becomes a lattice. The general idea of the algorithm is to compute 
 $\Pr(u \in H_{\alpha})$, sample $\alpha \in {\Bbb Z}$ from the computed probability distribution, assign $\xi_n=\alpha$ and then iterate, until all coordinates are sampled. 
 
 We will use the following inequality from \cite{Ba03} and \cite{AR05}:
 $$\Theta_{\Lambda}(v) \ \leq \ \Theta_{\Lambda}(0) \ \leq \ \exp\left\{\dist^2(v, \Lambda)\right\}
  \Theta_{\Lambda}(v) \quad \text{for all} \quad v \in {\Bbb R}^n. \tag7.1.3$$

 The following lemma summarizes various technical estimates that we need.
  
 \proclaim{(7.2) Lemma} Let $v_{\alpha}$ be the orthogonal projection of $v$ onto $H_{\alpha}$, so that 
 $$\|v-v_{\alpha}\|=\dist(v, H_{\alpha}).$$
 \roster
 \item We have
$$\Pr(\xi_n=\alpha)=\Pr(u \in H_{\alpha}) = \exp\left\{ - \|v-v_{\alpha}\|^2\right\} {\Theta_{\Lambda_{\alpha}}(v_{\alpha}) \over \Theta_{\Lambda}(v)};$$
\item We have 
$$\Pr(\xi_n=\alpha) \ \leq \ \exp\left\{ \dist^2(v, \Lambda) -\|v-v_{\alpha}\|^2 \right\};$$
\item Let $B$ be the Gram matrix of $u_1, \ldots, u_n$ and suppose that 
$$\lambda_{\min} I \ \preceq \ B \ \preceq \ \lambda_{\max} I$$
for some $ \lambda_{\max} \geq \lambda_{\min} >0$. Let
$$v=\eta_1 u_1 + \ldots + \eta_n u_n$$
for some real $\eta_1, \ldots, \eta_n$. Then 
$$\Pr(\xi_n=\alpha)\ \leq \ \exp\left\{ {n \lambda_{\max}  \over 4} - \lambda_{\min} (\eta_n - \alpha)^2\right\}.$$
\item Suppose that the Gram matrix $B$ of $u_1, \ldots, u_n$ satisfies the condition of Part (3) for some $\lambda_{\max} \geq \lambda_{\min} >0$.
 Then the Gram matrix $B'$ of $u_1, \ldots, u_{n-1}$ satisfies the condition with the same $\lambda_{\max}$ and $\lambda_{\min}$.
 \endroster
 \endproclaim
 \demo{Proof} For every $u \in H_{\alpha}$, by the Pythagoras Theorem, we have 
 $$\|u-v\|^2 = \|v- v_{\alpha}\|^2 + \|v_{\alpha} - u\|^2.$$
 Hence
 $$\sum_{u \in \Lambda_{\alpha}} \exp\left\{-\|u-v\|^2 \right\} = \exp\left\{-\|v-v_{\alpha}\|^2\right\} \sum_{u \in \Lambda_{\alpha}} \exp\left\{ - \|u-v_{\alpha}\|^2 \right\},$$
 and the proof of Part (1) follows.
 
 To prove Part (2), by applying (7.1.3) we get 
 $$\Theta_{\Lambda_{\alpha}}\left(v_{\alpha}\right) \ \leq \ \Theta_{\Lambda_{\alpha}}(0)=\Theta_{\Lambda_0}( 0) \ \leq \ \Theta_{\Lambda}( 0) 
 \ \leq \ \Theta_{\Lambda}(v) \exp\left\{ \dist^2(v, \Lambda)\right\},$$
 and the proof follows from Part (1).
 
 Next, we prove Part (3). For $i=1, \ldots, n$, let $\nu_i$ be the integer nearest to $\eta_i$, so that $|\eta_i - \nu_i | \leq {1 \over 2}$ and let 
 $u=\nu_1 u_1 + \ldots + \nu_n u_n$, so that $u \in \Lambda$. 
 Let 
 $$y=\left( \eta_1 - \nu_1, \ldots, \eta_n -\nu_n \right).$$
  Then 
 $$\dist^2(v, \Lambda) \ \leq \ \|v-u\|^2 = \langle By, y \rangle \ \leq \ \lambda_{\max} \|y\|^2 \ \leq \ {\lambda_{\max} n \over 4}. \tag7.2.1$$
 Let $w$ be a unit vector orthogonal to $u_1, \ldots, u_{n-1}$. Then 
 $$\|v-v_{\alpha}\|^2=\dist^2(v, H_{\alpha})=\left( \langle v, w \rangle - \langle v_{\alpha}, w \rangle \right)^2 =\langle u_n, w \rangle^2 (\eta_n - \alpha)^2. \tag7.2.2$$
 To bound $\langle u_n, w \rangle^2$, we consider the $n \times n$ matrix $A$ having vectors $u_1, \ldots, u_n$ as rows. Then $B=A A^T$ and since the eigenvalues of the matrices $A A^T$ and $A^T A$ coincide (the matrices are similar), we also have 
 $$\lambda_{\min} I \ \preceq \ A^T A, \tag7.2.3$$
 Now, $A w=\langle u_n, w \rangle e_n$, where 
 $e_n$ is the $n$-th standard basis vector and hence $A^T A w = \langle u_n, w \rangle u_n$. From (7.2.3), we obtain that
 $$\langle A^T Aw, w \rangle =\langle u_n, w \rangle^2 \ \geq \ \lambda_{\min}. \tag7.2.4$$
Combining (7.2.1), (7.2.2), (7.2.4) and Part (2), we complete the proof of Part (3).

To prove Part (4), we identify ${\Bbb R}^{n-1}$ with the coordinate subspace of ${\Bbb R}^n$, consisting of the points $x=\left(\xi_1, \ldots, \xi_n\right)$ where 
$\xi_n=0$. The condition on the matrix $B$ says that 
$$\lambda_{\min} \|x\|^2 \ \leq \ \langle Bx, x \rangle \ \leq \ \lambda_{\max} \|x\|^2 \quad \text{for} \quad x \in {\Bbb R}^n,$$ 
while the same condition for $B'$  says that the above inequality holds for $x \in {\Bbb R}^{n-1} \subset {\Bbb R}^n$.
 {\hfill \hfill \hfill} \qed
 \enddemo
 
 Now we are ready to present the sampling algorithm.
 
 \subhead (7.3) Algorithm for sampling from the discrete Gaussian distribution \endsubhead
 \bigskip
 {\bf Input:} A basis $u_1, \ldots, u_n$ of a lattice $\Lambda$ such that the Gram matrix $B$ of $u_1, \ldots, u_n$ satisfies 
 $$\lambda_{\min} I \ \preceq \ B \ \preceq \ \lambda_{\max} I$$
 for some $\lambda_{\max} \geq \lambda_{\min} > 0$ such that 
 $$\lambda_{\min} \ \geq \ \pi^2\left(s+ {e^s \over 4}\left(1-e^{-s}\right)\left(1-e^{-2s}\right)\right)^{-1}\quad \text{and} \quad 
 \lambda_{\max} \ \leq \ \pi^2 s^{-1}$$
 for some $s \geq 1$, a point $v \in {\Bbb R}^n$ and $0 < \epsilon \leq 1$.
 \medskip
 {\bf Output:} A random point $u$ from a distribution $\mu$ on $\Lambda$ such that 
 $${1 \over 2} \sum_{u \in \Lambda} \left| \mu(u) -\Pr(u)\right| \ \leq \ \epsilon, \quad \text{where} \quad \Pr(u)=
 {\exp\left\{-\|u-v\|^2\right\} \over \Theta_{\Lambda}(v)}.$$
 \medskip
 {\bf Algorithm:} 
 \smallskip
 {\bf Step 0:} Let 
 $$v=\eta_1 u_1 + \ldots + \eta_n u_n.$$
 From Part (3) of Lemma 7.2, compute an integer $l \geq 1$, 
 $$l=O\left({n \over \lambda_{\min}} \ln {n \over \epsilon}\right),$$ such that for $u =\xi_1 u_1 + \ldots + \xi_n u_n$, $u \in \Lambda$, one has
$$\Pr\Bigl( |\xi_i -\eta_i| \ > \ l \quad \text{for some} \quad i=1, \ldots, n \Bigr) \ < \ {\epsilon \over 10n}.$$
 For $k=1, \ldots, n$ perform the following steps.
 \smallskip
 {\bf Step $k$:} The input of Step $k$ is the lattice $\Lambda^{(k)} \subset {\Bbb R}^{n-k+1}$ with basis \newline $u_1, \ldots, u_{n-k+1}$, where 
 ${\Bbb R}^{n-k+1}$ is identified with $\spa\left(u_1, \ldots, u_{n-k+1}\right)$, and a point $v^{(k)} \in {\Bbb R}^{n-k+1}$,
 $$v^{(k)}=\eta_1^{(k)} u_1^{(k)} + \ldots + \eta_{n-k+1}^{(k)} u_{n-k+1}.$$
 When $k=1$, we have 
 $\Lambda^{(1)}=\Lambda$ and $v^{(1)}=v$. For $\alpha \in {\Bbb Z}$ such that 
 $$| \alpha - \eta_{n-k+1}| \ \leq \ l,$$
 compute the probabilities that $\xi_{n-k+1}=\alpha$ within relative error $\epsilon/10n$ as in Part (1) of Lemma 7.2. To compute theta functions, use the algorithm of Section 3.2 and the reciprocity relation (1.1.4).
 Sample a value
 $\xi_{n-k+1}=\alpha$ from the resulting probability distribution. If $k< n$, let $v^{(k+1)} = v_{\alpha}^{(k)}$ and go to Step $k+1$.
 \smallskip
 At the end of Step $n$, we have integers $\xi_1, \ldots, \xi_n$. Output 
 $$u =\xi_1 u_1 + \ldots + \xi_n u_n.$$
 \bigskip
 We state the result as a theorem.
 \proclaim{(7.4) Theorem}  The algorithm of Section 7.3 samples a point $u \in \Lambda$ from a distribution which is $\epsilon$-close in the total variation distance to the discrete Gaussian distribution (7.1.1) in time polynomial in $n$, $\epsilon^{-1}$ and $\lambda_{\min}^{-1}$.
 \endproclaim
 {\hfill \hfill \hfill} \qed 
 
 \remark{(7.5) The smooth case} As we mentioned in Section 2.4, the algorithm follows the general scheme of Peikert \cite{Pe10}. The difference is that \cite{Pe10} deals with the smooth range, when $B \preceq sI$ with $s\ll (\ln n)^{-1}$ so that the value of $\Theta_{\Lambda}(v)$ does not significantly depend on the choice of $v\in {\Bbb R}^n$. Hence there is no need to compute values of the theta function, and one needs to sample $\alpha$ from 
 the distribution where 
 $$\Pr\left(\xi_{n-k+1}=\alpha\right) \sim \exp\left\{ - \|v^{(k)}-v^{(k)}_{\alpha}\|^2 \right\}. \tag7.5.1$$
 Another computational advantage of the smooth case in that the distribution (7.5.1) is well-approximated by a continuous Gaussian distribution. As a result, the complexity of sampling $\xi_{n-k+1}$ does not depend badly on the length of an interval for $\xi_{n-k+1}$ and so there is no dependence on 
 $\lambda_{\min}$ that we have in Theorem 7.4. It appears that once we leave the smooth range, we do need to compute theta functions, and the dependence on $\lambda_{\min}$ appears to be unavoidable.
  \endremark

\head 8. The smooth range \endhead

Let us fix $\gamma >1$.
In this section, we present a fully polynomial time approximation scheme (FPTAS) for computing (1.1.3) when $B$ is an $n \times n$ positive definite matrix of a sufficiently large size $n \geq n_0(\gamma)$ satisfying
$$s I \ \preceq \ B \quad \text{where} \quad s \ \geq \ \gamma \ln n.$$
Thus we present a deterministic algorithm that for any $0 < \epsilon \leq 1$ approximates (1.1.3) within relative error $\epsilon$ in time polynomial in $\epsilon^{-1}$ and $n$. From the reciprocity relation (1.1.4), we immediately get an FPTAS for approximating $\Theta(B, y)$ provided 
$$B \ \preceq \ s I \quad \text{where} \quad s \ \leq \ {\pi^2 \over \gamma \ln n} I$$
as long as $n \geq n_0(\gamma)$.
The results of this section are likely to be known in some form, but since we are unable to provide a reference, we summarize them here for completeness.

The algorithm is based on the following simple result.

\proclaim{(8.1) Theorem} Fix $\gamma >1$ and 
let $B$ be an $n \times n$ be positive definite matrix such that 
$$sI \ \preceq \ B \quad \text{where} \quad s \ \geq \ \gamma \ln n.$$
\roster 
\item 
For $n \geq 2$ and for all integer $k \geq 1$, we have 
$$\sum\Sb x \in {\Bbb Z}^n: \\ \|x\|^2 \geq  k \endSb \exp\left\{ - \langle B x, x \rangle \right\} \ \leq \ 60 n^{(1-\gamma)k}.$$
\item Let 
$$n_0(\gamma)= \exp\left\{{5 \over \gamma-1}\right\}.$$
Then for any $n \geq n_0(\gamma)$
and any $b \in {\Bbb R}^n$, we have 
$$\left| -1 + \sum_{x\in {\Bbb Z}^n} \exp\left\{ -\langle Bx, x \rangle + \ii \langle b, x \rangle\right\} \right| \ \leq \ {1 \over 2}.$$
\item For any integer $k \geq 1$, we have 
$$\left|x \in {\Bbb Z}^n: \ \|x\|^2 \leq k \ \right| \ \leq \ (2n+2)^k.$$
\endroster
\endproclaim
\demo{Proof} The proof of Part (1) is similar to that of Lemma 5.4. For $0 < \tau < s$, we have
$$\split &\sum\Sb x \in {\Bbb Z}^n: \\ \|x\|^2 \geq k \endSb \exp\left\{ -\langle Bx, x \rangle \right\} \ \leq \ \sum\Sb x \in {\Bbb Z}^n: \\ \|x\|^2 \geq k \endSb e^{-s \|x\|^2} \ \leq \ e^{-\tau k} \sum\Sb x \in {\Bbb Z}^n: \\ \|x\|^2 \geq k \endSb  
e^{-s\|x\|^2} e^{\tau \|x\|^2} \\ &\qquad \leq \ e^{-\tau k} \Theta\bigl((s-\tau)I\bigr) \ \leq \ \exp\left\{ - \tau k + {2ne^{-(s-\tau)} \over 1-e^{-(s-\tau)}}\right\},
\endsplit$$
where the last inequality is from Part (1)  of Lemma 5.4. We choose 
$$\tau = (\gamma -1) \ln n.$$ 
Since $s-\tau \geq \ln n$, we obtain
$$ \exp\left\{ - \tau k + {2ne^{-(s-\tau)} \over 1-e^{-(s-\tau)}}\right\} \leq  \exp\left\{ -\tau k + {2 \over 1- n^{-1}}\right\} \ \leq \ 60 n^{(1-\gamma) k},$$
which completes the proof of Part (1).

Part (2) follows from Part (1) for $k=1$, since for $n \geq n_0(\gamma)$ we have $60 n^{1-\gamma} \leq {1 \over 2}$.

To prove Part (3), letting $x=\left(\xi_1, \ldots, \xi_n\right)$ and $\eta_i=\xi_i^2$, we observe that the number  non-negative integer solutions to the inequality 
$\eta_1 + \ldots + \eta_n \leq k$ is 
$${n+k \choose k}={(n+k)(n+k-1) \cdots (n+1) \over k(k-1) \cdots 1} \ \leq \ (n+1)^k.$$
Since each of at most $k$ positive $\eta_i$ correspond to at most two values $\pm \xi_i$, the bound follows.

{\hfill \hfill \hfill} \qed
\enddemo

Now we are ready to present the algorithm.

\subhead (8.2) The algorithm \endsubhead Fix $\gamma >1$ and 
$$n_0=\exp\left\{ 5 \over \gamma-1\right\}.$$
\bigskip
{\bf Input:} For $n \geq n_0(\gamma)$, an $n \times n$ positive definite matrix $B$ such that $s I \preceq B$ for some $s \geq \gamma \ln$,
a vector $b \in {\Bbb R}^n$ and $0 < \epsilon < 1$.
\medskip
{\bf Output:} A number approximating 
$$\sum_{x \in {\Bbb Z}^n} \exp\left\{ - \langle Bx, x \rangle +   \ii \langle b,  x \rangle\right\} \tag8.2.1$$
within relative error $\epsilon$.
\medskip 
{\bf Algorithm:} From Parts (1) and (2) of Theorem 1, choose 
$$k=O\left({\ln (1/\epsilon) \over (\gamma-1) \ln n} \right),$$
so that
$$\sum\Sb x \in {\Bbb Z}^n: \\ \|x\|^2 \leq k\endSb \exp\left\{ -\langle Bx, x \rangle + \ii \langle b, x \rangle \right\} \tag8.2.2$$
approximates (8.2.1) within relative error $\epsilon$, and compute (8.2.2).

From Part (3) of Theorem 8.1, the sum (8.2.2) contains $(1/\epsilon)^{O\left({1\over \gamma-1}\right)}$
 terms.

\head Acknowledgment \endhead

The author is grateful to the anonymous referees for their criticism and suggestions.

\enddocument